 \magnification=1200
\input amssym.def
\input amssym.tex
 \font\newrm =cmr10 at 24pt
\def\bul{\raise .9pt\hbox{\newrm .\kern-.105em } }

 \def\fr{\frak}
 \font\sevenrm=cmr7 at 7pt

 \baselineskip 20pt
 
 \def\h{\hbox{ }}

 \def\u{{\fr u}}
 \def\r{{\fr r}}

 \def\n{{\fr n}}

 \def\ss{{\fr s}}
 
 \def\b{{\fr b}}
 
 \def\hh{{\fr h}}
 
 \def\ee{{\fr e}}

 \def\g{{\fr g}}

 \def\e{{\tilde e}}

 \def\<{\le}
 \def\>{\ge}

 \def\s{{\h\subset\h}}
 
 \def\vs{\vskip }

 \def\mapright#1
  {\smash{\mathop
  {\longrightarrow}
  \limits^{#1}}}

 \def\kk#1{{\kern .4 em} #1}
 \def\vs{\vskip 1pc}

\hsize = 30pc
\vsize = 45pc
\overfullrule = 0pt

\rm

	\centerline{\bf The cascade of orthogonal roots and the coadjoint
structure of the}\centerline {\bf nilradical of a Borel
subgroup of a semisimple Lie group}\vskip 1.5pc 
\centerline{\it To the memory of a dear friend, I. M. Gelfand, one
of the}
\centerline{\it great mathematicians of the 20th century}\vskip
1.5pc
\centerline{Bertram Kostant}\vskip 1.5pc 

{\bf Abstract:} \rm Let $G$ be a semisimple Lie group and let $\g =
\n_- +
\hh +\n$ be a triangular decomposition of $\g= \hbox{Lie}\,G$. Let
$\b =
\hh +\n$ and let $H,N,B$ be Lie subgroups of $G$ corresponding
respectively to $\hh,\n$ and $\b$. We may identify $\n_-$ with the
dual space to $\n$. The coadjoint action of $N$ on $\n_-$ extends
to an action of $B$ on $\n_-$. There exists a unique nonempty
Zariski open orbit $X$ of $B$ on $\n_-$. Any $N$-orbit in $X$ is a
maximal coadjoint orbit of $N$ in $\n_-$. The cascade of orthogonal
roots defines a cross-section $\r_-^{\times}$ of the set of such
orbits leading to a decomposition $$X = N/R\times \r_-^{\times}.$$
This decomposition, among other things, establishes the structure of
$S(\n)^{\n}$ as a polynomial ring generated by the prime
polynomials of $H$-weight vectors in $S(\n)^{\n}$. It also leads to
the multiplicity 1 of $H$ weights in $S(\n)^{\n}$. \vskip 6pt

\vskip 6pt
{\bf MSC (2010) subject codes:} representation theory, invariant theory.
\vskip 6pt
{\centerline {\bf 0. Introduction} {\bf 0.1}. Let 
$\g$ be a complex semisimple Lie algebra and let
$\hh$ be a Cartan subalgebra of $\g$. The Killing form $(x,y)$, 
 denoted by ${\cal K}$, on
$\g$ induces a nonsingular bilinear form $(\mu,\nu)$ on the dual
space
$\hh^*$ to $\hh$. Let $\Delta\s \hh^*$ be the set of roots
corresponding to $(\hh,\g)$. For each $\varphi\in\Delta$ let
$e_{\varphi}\in\g$ be a corresponding root vector. Let $\b$ be a Borel
subalgebra of $\g$ which contains $\hh$ and let $\n$ be the nilradical of $\b$. Let $$\n_- + \hh +\n$$
be a triangular decomposition of $\g$. Then a choice $\Delta_+$ (respectively $\Delta_-$) of positive
resp. negative) roots is chosen so that $\Delta_+ = \{\varphi\mid e_{\varphi}\in \n\} $ (resp.
$\{\varphi\mid e_{-\varphi}\in \n_-\} $). One has $\Delta_- = -\Delta_+$. If $\ss\s \g $ is a Lie
subalgebra, then $S(\ss)$ and $U(\ss)$ will denote respectively the symmetric and
enveloping algebras of $\ss$. We are mainly concerned with the case where $\ss = \n$
and with the structure of the space of
$\n$-invariants $S(\n)^{\n}$ (or equivalently $\hbox{cent}\,U(\n)$) and with the
action of $\hh$ on
$S(\n)^{\n}$ and $\hbox{cent}\,U(\n)$. 
 
Let $G$ be a Lie group such that $\g = \hbox{Lie}\,G$ and let $H$, $N$, and $B$ be the
Lie subgroups which correspond respectively to $\hh$,$\n$ and $\b$. Using ${\cal K}$
we may identitfy
$\n_-$ as the dual space to $\n$ so that the coadjoint action of $N$ defines an action of $N$ on $\n_-$.
Furthermore since $B$ normalizes $N$ there is accompanied an action of $B$ (and in a particular $H$) on
$\n_-$. Many of the results in this paper are quite old and were cited in [J]. The techniques in this
paper are algebra-geometric in nature and are quite different from those in $[J]$. 
We will see that
$B$ has a unique Zariski open and (dense) orbit $X$ in $\n_-$. The results arise from the rather
elegant structure of the coadjoint action of $N$ on the affine variety $X$ and from the equally elegant
action of
$H$ on these
$N$-coadjoint orbits. The main tool is the use of our well-known cascade of
orthogonal roots. The cascade plays a major role in a number of papers. In particular
the cascade has been used in [J] and [L-W].
We will recall the definition of the cascade ${\cal B}$ in
\S 1 below and elaborate on some of its properties. Among other things ${\cal B}$ is a special
maximal set of strongly orthogonal roots. Let $m = \hbox{card}\,{\cal B}$ and we will
write ${\cal B} = 
\{\beta_1,\ldots,\beta_{m}\}$.

 Let $\r_-\s \n_-$ be the span of $\{e_{-\beta_i}\},\,i=1,\ldots,m$, and let $\r_-^{\times}$ be the Zariski
open subset of $\r_-$ defined by the condition that if $e\in \r_-$, then $e\in\r_-^{\times}$ if all the
coefficients of $e$ relative to the $e_{-\beta_i}$ are nonzero. Let $T$ be the toroidal subgroup of $H$
defined so that $\beta_i^{\vee}, i=1,\ldots,m$, is a basis of $\hbox{Lie}\,T$. Then
$T$ operates on $\r_-$ and
$\r_-^{\times}$ is the unique Zariski open orbit of $T$ in $\r_-$. It is clear that $T$ operates simply
and transitively on $\r_-^{\times}$ so that as affine varieties $T\cong \r_-^{\times}$. In this paper we will
prove {\bf (1)} that $\r_-^{\times} \s X$ and {\bf (2)} $\r_-^{\times}$ is a
cross-section of the set of $N$-coadjoint orbits in $X$. In addition, {\bf (3)} every
$N$ orbit on $X$ is a maximal coadjoint orbit of $N$ in $\n_-$ and furthermore {\bf(4)}, the isotropy
group $R$ of
$N$ at any point $p$ in $\r_-^{\times}$ is independent of $p$. Moreover if $A(Y)$
denotes the affine ring of an affine variety
$Y$, then {\bf (5)} one has, as affine varieties, $X\cong N/R\times \n_-^{\times}$ so that $$A(X) \cong
A(N/R)\otimes A(T)\eqno (0.1)$$ and {\bf (6)}
$$A(X)^N
\cong A(T),\eqno (0.2)$$ noting {\bf (7)}, that $A(T)$ is the character ring of the torus $T$ and is
generated by the cascade
${\cal B}$ which, in fact, is a basis of the character ring. 

Returning to the description of $S(\n)^{\n}$ we establish that {\bf (a)} $S(\n)^{\n}$ embeds naturally in
$A(T)$ and that $S(\n)^{\n}$ and $A(T)$ have the same quotient field. Let $\widehat {T}$ be the
character group of $T$ so that $\widehat {T}\s A(T)$. For any $\xi\in \widehat {T}$ let $\nu(\xi)$
be the corresponding $H$-weight so that $$\{\nu(\xi)\mid \xi\in \widehat {T}\}\,\,\hbox{ is the
free abelian group generated by ${\cal B}$}. $$ Let
$Q= S(\n)^{\n}\cap \widehat {T}$. One readily has that $\nu(\xi)$ is dominant for any $\xi\in Q$. Now
let $P$ be the set of all $\xi\in Q$ which, as polynomials on $\n_-$, are prime. We prove {\bf (b)}
$\hbox{card}\,P = m$ and {\bf (c)} if $P= \{\xi_1,\ldots,\xi_m\}$,  then for $i=1,\ldots,m$, the $\xi_i$
are algebraically independent and $\nu_i = \nu(\xi_i)$ are linearly independent. Furthermore if $d =
(d_1,\ldots,d_m)\in
\Bbb Z^{m}$ and $$\xi^d = \xi_1^{d_1}\cdots \xi_m^{d_m}, \eqno (0.3)$$ note that $$\nu(\xi^d) =
\sum_{i=1}^m d_i\,\nu_i.\eqno (0.4)$$ We prove {\bf (d)}
$$\widehat {T} =\{\xi^d\mid d\in \Bbb Z^m\},\eqno (0.5)$$ and if $\Bbb N = \{0,1\,\ldots,\}$ is the set
of natural numbers, then we also establish {\bf (e)}
$$Q =
\{\xi^d\mid d\in \Bbb N^m\}\eqno (0.6)$$ so that {\bf (f)} $$S(\n)^{\n}\,\,\hbox{is the polynomial ring
$\Bbb C[\xi_1,\ldots,\xi_m]$ and $Q$ is a weight basis of this ring}.\eqno (0.7)$$ In particular
{\bf (g)} every weight in $Q$ has multiplicity 1 and if $\eta\in Q$, then {\bf (h)} writing $$\eta=
\xi^d,\,\,\hbox{where
$d\in \Bbb N^m$, is the prime decomposition of the polynomial $\eta$}. \eqno (0.8)$$ Furthermore since any
$\nu(\eta)$ is dominant for any $\eta\in Q$ the coefficients $k_j$ in the expansion $\nu(\eta)=
\sum_{j=1}^m k_j\beta_j$ are nonnegative integers and {\bf (i)} as a polynomial on $\n_-$,
  $$\hbox{deg}\,\eta =
\sum_{j=1}^m k_j.\eqno (0.9)$$ \vskip .2pc 

{\bf Remark 0.1.} Given a dominant weight, we constructed in [K]
(modifying the method of Lipsman--Wolf) an element
$f_{(k)}\in S(\n)^{\n}$ of degree $k$. If the dominant weight is $\rho$, equal to one-half the sum of 
the positive roots, then $\nu(f_{(k)}) = 2\rho$, and one readily shows that all the elements $\xi_i$ in
$P$  may be given as the prime factors of $f_{(k)}$. \vskip .2pc 

Finally returning to (0.1) one establishes {\bf (j)} that
$S(n)^{\n}$ ``separates" all the maximal
$N$-coadjoint orbits that lie in $X$. \vskip .5pc

 {\bf Remark 0.2.} A number of results in [J] were unknown
to us when
$[J]$ was written. For example we were unaware that the set $$\{\nu(\eta)\mid \eta\in Q\}$$ included all the
dominant elements in the lattice $L$ generated by ${\cal B}$. This, however, follows from an argument in [K]
which proves that if
$\eta\in L$ and $\nu(\eta)$ is dominant, then $\eta^4\in Q$, which by (0.6), implies that $\eta\in Q$.
Joseph in [J] does not deal with prime polynomials. Instead he sets up a bijection of ${\cal B}$ with
certain generators  of
$S(\n)^{\n}$, $$\beta_j\mapsto \eta_j.\eqno (0.10)$$ 
Furthermore he very cleverly determines, by a sort of
Gram--Schmidt process, the $m\times m$ matrix $s_{ij}$ where $$\nu(\eta_j) = \sum_{i=1}^m s_{ij}\beta_i\eqno
(0.11)$$ In addition very useful information is given in his tables II and III. Among other things
the $\nu(\eta_j)$ are expressed in terms of the fundamental representations of $\g$.\vskip 6pt {\bf 0.2.}
The results in this paper were inspired by Dixmier's result [D] for the special case where $G = Sl(n,\Bbb
C)$. \vskip 1pc
\centerline{\bf 1. The cascade of orthogonal roots}\vskip 2pt 
  {\bf 1.1.} Let $\ell = \hbox{rank}\,\g$ and let $\Pi\s \Delta_+$ be the set of
simple positive roots. For each $\varphi\in \Delta_+$ let
$n_{\alpha}(\varphi)\in \Bbb Z_+$ be such that
$\varphi=\sum_{\alpha\in\Pi}n_{\alpha}(\varphi)\,\alpha$. Now let
$\Pi(\varphi)=\{\alpha\in \Pi\mid n_{\alpha}(\varphi)>0\}$.
Then, as one knows, $\Pi(\varphi)$ is a connected subset of $\Pi$.
Hence there is a unique complex simple Lie subalgebra $\g(\varphi)$
of
$\g$, with Cartan subalgebra $\hh$, having $\Pi(\varphi)$ as a set
of simple positive roots. Let $\Delta(\varphi)\s \Delta$ be the 
set of roots of $(\hh,\g(\beta))$ and $\Delta(\varphi)_+ =
\Delta(\varphi)\cap \Delta_+$. Let $\b(\varphi) =
\b\cap\g(\varphi)$ and $\n(\varphi) = \n\cap \g(\varphi)$.

Let $\beta\in
\Delta_+$. We will say that $\beta$ is locally high if $\beta$ is
the highest root of $\g(\beta)$. If $\beta\in \Delta_+$ is locally
high,  let $E(\beta) = \{\varphi\in \Delta(\beta)\mid (\varphi,\beta)>0\}$ and let $\ee(\beta)$ be the span of
$e_{\varphi}$ for $\varphi\in E(\beta)$. Let $h^{\vee}(\beta)$
be the dual Coxeter number of $\g(\beta)$. Then, regarding the one-dimensional
  Lie algebra as a Heisenberg Lie
algebra, one knows 
\vskip 3pt {\bf Proposition 1.1.} {\it
$$\eqalign{&(1)\,\,\ee(\beta)\s\n(\beta)\,\,
\hbox{and $\ee(\beta)$ is an ideal in $\b(\beta)$}\cr
&(2)\,\,\ee(\beta)\,\,\hbox{is a Heisenberg Lie algebra of
dimension
$2h^{\vee}(\beta)-3$ and with center $\Bbb C
e_{\beta}$}\cr &(3)\,\,2(\beta,\varphi)/(\beta,\beta)
=1,\,\,\forall
\varphi\in E(\beta)/\{\beta\}.\cr}$$}\vskip 2pt For any $\varphi\in
\Delta$ let $\u(\varphi)$ be the TDS spanned by
$h_{\varphi},e_{\varphi}$ and $e_{-\varphi}$ and let
$\g(\beta)^o$ be the semisimple Lie subalgebra (possibly zero) of
$\g(\beta)$ spanned by all $\u(\varphi)$ where $\varphi\in
\Delta(\beta)$ and
$(\beta,\varphi)=0$. \vskip 6pt {\bf Remark 1.2.} One notes that the
highest root of any simple component of $g(\beta)^o$ is locally
high. It is necessarily  orthogonal to $\beta$. \vskip 3pt
Introduce the usual partial ordering in the weight
lattice where $$\nu> \mu \eqno (1.1)$$ if $\nu-\mu$ is a sum of
positive roots.

A sequence $$C=
\{\beta_1,\ldots,\beta_k\}\eqno (1.2)$$ of positive roots will be
called a chain cascade if $\beta_1$ is the highest root of a simple
component of
$\g$ and, inductively, if
$1<j\leq k$, and
$\beta_i$ has been given for $1\leq i<j$ and is locally high
for all such $i$ then
$\beta_j$ is the highest root of a simple component of
$\g(\beta_{j-1})^o$. \vskip 6pt {\bf Remark 1.3.} One notes a chain cascade is simply ordered.
For the chain cascade above $\beta_1$ is the maximal element and $\beta_k$ is the minimal element.
\vs Let ${\cal B}$ be the set (cascade) of all positive roots
$\beta$ which are members of some chain cascade. As an immeditate
consequence of Remark 1.2 one has
\vskip 6pt {\bf Proposition 1.4.} {\it Any root $\beta$ in ${\cal B}$ is locally
high.}\vs Clearly for any positive root $\varphi$ one uniquely
defines a chain cascade $$C(\varphi)=
\{\beta_1,\ldots,\beta_k\}\eqno (1.3)$$ by the condition that (1)
$\varphi\in
\Delta(\beta_i)_+,\,i=1,\ldots,k$, and (2)
$(\varphi,\beta_i)=0,\,i=1,\ldots,k-1$, and $(\varphi,\beta_k) >0$. 
In the notation of (1.2) clearly $$C = C(\beta_k).$$ One
readily also notes
\vskip 6pt {\bf Proposition 1.5.} {\it  For
$\varphi,\varphi'\in \Delta_+$ one has $$C(\varphi) =
C(\varphi')\,\,\iff\,\,\hbox{there exists $\beta\in {\cal B}$ such that
$\varphi,\varphi'\in E(\beta)$}\eqno (1.4)$$ in which case $\beta$ is the
minimal element of $C(\varphi) =
C(\varphi')$. In particular one has the disjoint union $$\Delta_+ =
\cup_{\beta\in {\cal B}}\,\,E(\beta)\eqno (1.5)$$ and the consequential
direct sum (as a vector space) of Heisenbergs
$$\n= \oplus_{\beta\in {\cal B}}\,\, \ee(\beta).\eqno (1.6)$$}\vskip 6pt {\bf
Lemma 1.6.} {\it Any two distinct elements of a chain cascade $C$
are strongly orthogonal.} \vskip 6pt {\bf Proof.} By definition it is clear
that any distinct members of $C$ are orthogonal. Without loss assume $C$ is
given by (1.2) and $1\leq i<j\leq k$. Then $\beta_j\in
\Delta(\beta_i)$. But then $\beta_i+ \beta_j$ cannot be a root since
$\beta_i$ is the highest root of $\g(\beta_i)$. \hfill QED
\vskip 6pt Given
a chain cascade $C$, say given by (1.2), it is clear that $C' =
\{\beta_1,\ldots,\beta_j\}$ is a chain cascade for any $j\leq k$.
Given two chain cascades $C$ and $C'$ we will say that $C'$ is a
subchain of $C$ if $C$ and $C'$ are of this form. Also subsets $\Delta^1$ and $\Delta^2$ of $\Delta$ will
  be called totally disjoint if any element of $\Delta_1$ is strongly
orthogonal to every element of $\Delta_2$. \vskip 6pt {\bf Proposition
1.7.} {\it Let $\beta,\beta'\in {\cal B}$ be distinct. If $C(\beta)$
is not a subchain of $C(\beta')$ and vice versa,  then not
only are
$\beta$ and $\beta'$ strongly orthogonal but in fact $\Delta(\beta)$
and $\Delta(\beta')$ are totally disjoint.} \vskip 6pt {\bf Proof.}
Without loss assume that $\beta = \beta_k$ in the notation of (1.2)
so that in that notation $C(\beta) = C$. Let $C(\beta') =
\{\beta'_1,\ldots,\beta'_{k'}\}$. Without loss we may assume that
$k\leq k'$. By our asumption on subchains one has
$\beta'_k\neq \beta_k$. Let $j\leq k$ be minimal such that
$\beta'_j\neq \beta_j$. If $j=1$,  then the result is clear since
$e_{\beta}$ and $\e_{\beta'}$ lie in different simple components of
$g$. Assume $j>1$. Then $\beta_{j-1} = \beta'_{j-1}$, but then
$e_{\beta}$ and $e_{\beta'}$ lie in different simple components
of $\g(\beta_{j-1})^o$. The result then follows.\hfill QED
\vskip 6pt {\bf Theorem
1.8.} {\it The set ${\cal B}$ is a maximal set of strongly orthogonal
roots.} \vskip 6pt {\bf Proof.} That ${\cal B}$ is a set of strongly orthogonal
roots follows from Lemma 1.6 and Proposition 1.7. That it is maximal
follows from the disjoint union (1.5).\hfill QED
\vskip 6pt We call
${\cal B}$ the cascade of strongly orthogonal roots.
\vskip 6pt {\bf 1.4.} Let $W$ be the Weyl group of $\g$ operating in $\hh$
and
$\hh^*$. For $\beta\in {\cal B}$ let $W(\beta)\s W$ be the Weyl group of
$\g(\beta)$. Reluctantly submitting to common usage, let $w_o$ be
the long element of $W$ and let $w_o(\beta)$ be the long element of
$W(\beta)$. For any $\beta\in {\cal B}$ let $s_{\beta}\in W$ be the
reflection defined by $\beta$. \vskip 6pt {\bf Proposition 1.9.} {\it Let
$\beta,\beta'\in {\cal B}$. Then $\Delta(\beta)$ is stable under
$s_{\beta'}$. Furthermore if $s_{\beta'}|\Delta(\beta)$ is not the
identity,  then $s_{\beta'}\in W(\beta)$.}

 \vskip 6pt {\bf Proof.} First
assume there exists a chain cascade $C$ containing both $\beta$ and
$\beta'$. Without loss assume $C$ is given by (1.2) and $\beta =
\beta_i$ and $\beta' = \beta_j$ for some $i,j\leq k$. 
Then $\Delta(\beta)$ is clearly element-wise fixed by $s_{\beta'}$
if $j<i$ and $s_{\beta'} \in W(\beta)$ if $j\geq i$. Thus it
remains to consider only the case when $C(\beta)$ is not a subchain
of $C(\beta')$ and vice versa. But then the result follows from
Proposition 1.7. \hfill QED

\vskip 6pt Since the elements in ${\cal B}$ are orthogonal
to one another the reflections $s_{\beta}$ evidently commute
with one another. The long element
$w_o$ of the Weyl group
$W$ is given in terms of the product of these commuting
reflections. Let $\Delta_- = -\Delta_+$ and $\Delta(\varphi)_- =
-\Delta(\varphi)_+$ for any $\varphi\in \Delta_+$. 

\vskip 6pt {\bf
Proposition 1.10.} {\it One has
$$w_o= \prod_{\beta\in {\cal B}}\,\,s_{\beta}\eqno (1.7)$$ noting that the
order of the product is immaterial because of commutativity.
Furthermore $w_o$ stabilizes $\g(\beta)$ for any $\beta\in {\cal B}$ and
$$w_o|\g(\beta) = w_o(\beta)|\g(\beta).\eqno (1.8)$$} {\bf
Proof.} Let $\kappa\in W$ be given by the right side of (1.7).
Clearly $$\kappa(\Delta(\beta)) =
\Delta(\beta) \eqno (1.9)$$ for any $\beta\in {\cal B}$ by
by Proposition 1.9. But also clearly $$\kappa(\beta) = -\beta\eqno
(1.10)$$ for any
$\beta\in {\cal B}$ so that $\kappa$ carries the highest root of
$\g(\beta)$ to the lowest root of $\g(\beta)$. But $E(\beta)\s
\Delta(\beta)_+$. Hence $\kappa(E(\beta)) \s \Delta(\beta)_-$. In
particular $$\kappa(E(\beta)) \s \Delta_-\eqno (1.11)$$ for any
$\beta\in {\cal B}$. But then $\kappa(\Delta_+) = \Delta_-$ by (1.5).
Thus $\kappa = w_o$. But then one has (1.8) by (1.9). \hfill QED\vskip
1pc \centerline{\bf 2. The coadjoint action}

\vskip 3pt {\bf
2.1.} Let
$G$ be a simply connected Lie group where $\g= \hbox{Lie}\,G$ and let $N\s G$ be the subgroup
corresponding to $\n$. Let $\n_-$ be the span of all
$e_{-\varphi}$ for $\varphi\in \Delta_+$. One has the direct sum
$$\g = \n_-\oplus \b.\eqno (2.1)$$ Let $\Phi:\g\to \n_-$ be the
projection with kernel $\b$. Using the Killing form
$(x,y)$ we identify the dual space $\n^*$ to $\n$ with $\n_-$. Let
$v\in \n,\,u \in N$ and $z\in \n_-$. The coadjoint action of
$v$, $\hbox{coad}\,v$, and $u$, $\hbox{Coad}\,u$, on $\n_-$ is given by
$$\eqalign{\hbox{coad}\,v(z)&= \Phi\,[v,z]\cr \hbox{Coad}\,u(z)&=
\Phi\,\hbox{Ad}\,u(z).\cr}\eqno (2.2)$$

Let $m = \hbox{card}\,{\cal B}$ and let $\r$ be the commutative $m$-dimensional
subalgebra of $\n$ spanned by $e_{\beta}$ for $\beta\in {\cal B}$. Let
$R\s N$ be the commutative unipotent subgroup corresponding to
$\r$. Let $\r_-\s \n_-$ be the span of $e_{-\beta}$ for $\beta\in
{\cal B}$. For any $z\in \r_-,\,\beta\in {\cal B}$, let $a_{\beta}(z)\in 
\Bbb C$
be defined so that $$z = \sum_{\beta\in
{\cal B}}\,a_{\beta}(z)\,e_{-\beta}\eqno (2.3)$$ and let $$\r_-^{\times} = 
\{\tau\in
\r_-\mid a_{\beta}(\tau)\neq 0,\,\,\forall
\beta\in {\cal B}\}.\eqno (2.4)$$ As an algebraic subvariety of $\n_-$
clearly 
$$\r_-^{\times} \cong \Bbb (C^{\times})^m.\eqno (2.5)$$

For any $z\in \n_-$ let $O_z$ be the $N$-coadjoint orbit
containing $z$. Let $N_z\s N$ be the coadjoint isotropy subgroup at
$z$ and let $\n_z = \hbox{Lie}\,N_z$. Since the action is algebraic
$N_z$ is connected and hence as $N$-spaces $$O_z\cong N/N_z.\eqno
(2.6)$$

Let $\ss\s \n$ be the span of $e_{\varphi}$ for $\varphi\in
\Delta_+/{\cal B}$ so that one has the vector space direct sum $$\n =
\r\oplus \ss\eqno (2.7)$$ and let $\ss_-$ be the span of
$e_{-\varphi}$ for $\varphi\in
\Delta_+/{\cal B}$ so that one has the vector space direct sum 
$$\n_- = r_-\oplus \ss_-.\eqno (2.8)$$ Also for any $\beta\in {\cal
B}$ let
$\ss(\beta)$ be the span of $e_{\varphi}$ for
$\varphi\in E(\beta)/\{\beta\}$ so that the Heisenberg $$\ee(\beta)
=
\ss(\beta)\oplus \Bbb C\,e_{\beta}.\eqno (2.9) $$ Let
$\ss(-\beta)$ be the span of $e_{-\varphi}$ for
$\varphi\in E(\beta)/\{\beta\}$. One has the
direct sums
$$\ss =
\oplus_{\beta\in {\cal B}}\,\ss(\beta)\eqno (2.10)$$ and $$\ss_-
=\oplus_{\beta\in {\cal B}}\,\ss(-\beta).\eqno (2.11)$$        
Let $\beta\in {\cal B}$. Since $\ee(\beta)$ is a Heisenberg Lie algebra,
given $\varphi\in E(\beta)/\{\beta\}$, there exists a unique
$\gamma\in E(\beta)/\{\beta\}$ such that $\varphi + \gamma =
\beta$. We refer to $\gamma$ as the Heisenberg twin to $\varphi$
and the pair $\{\varphi,\gamma\}$ as Heisenberg twins. The following
lemmas lead to a considerable simplification in dealing with the
coadjoint action of
$N$ on
$\n_-$.\vs {\bf Lemma 2.1.}. {\it Let $\beta\in {\cal B}$. Assume
$\varphi,\varphi'\in \Delta_+$ and $$\beta = \varphi +
\varphi'.\eqno (2.12)$$ Then $\varphi$ and $\varphi'$ are both in
$E(\beta)/\{\beta\}$ and are Heisenberg twins.} 

\vskip 6pt {\bf
Proof.} The equality (2.12) immediately implies that both $\varphi$
and $\varphi'$ are in $\Delta(\beta)$. But since
$(\beta,\beta)\neq 0$ the equality (2.12) also implies that $\varphi$
and $\varphi'$ cannot both be orthogonal to $\beta$. Hence at least
one of the two must be in $E(\beta)/\{\beta\}$. But then the other
is also in $E(\beta)/\{\beta\}$ by the existence of a Heisenberg
twin. \hfill QED

\vskip 6pt  {\bf Lemma
2.2.} {\it Let $\beta,\beta'\in {\cal B}$ and let $x\in \ss(\beta)$. Assume
$x\neq 0$. Then if $\beta'\neq \beta$, one has $$\hbox{coad}\,x(e_{-\beta'})
= 0.\eqno (2.13)$$ On the other hand if $y=\hbox{coad}\,x(e_{-\beta})$, then
$$y\neq 0\eqno (2.14)$$ and $$y\in \ss(-\beta).\eqno (2.15)$$}

\vskip 6pt   
{\bf Proof.} Let $\varphi\in E(\beta)/\{\beta\}$ and let 
$$z = \hbox{coad}\,e_{\varphi}(e_{-\beta'}).\eqno (2.16)$$ Assume $z\neq
0$. Then clearly there exists $\varphi'\in \Delta_+$ such that $z=
c\,e_{-\varphi'}$ for some
nonzero scalar $c$. But then  
$(e_{-\beta'},[e_{\varphi},e_{\varphi'}])\neq 0$. Hence $$\beta' =
\varphi + \varphi'.\eqno (2.17)$$ But then $\varphi$ and $\varphi'$
are Heisenberg twins in $E(\beta')/\{\beta'\}$ by Lemma 2.1. Thus
$z\neq 0$ implies $\beta = \beta'$. This proves (2.13) by writing 
$x$ as a sum of root vectors for roots in $E(\beta)/\{\beta\}$. 

Now assume $\beta = \beta'$. Let $\varphi'$ be the Heisenberg twin
to $\varphi$ in $E(\beta)/\{\beta\}$. But then clearly $z =
c\,e_{-\varphi'}$ for some
nonzero scalar $c$. But then (2.14) and (2.15) follow immediately.

\hfill QED
\vskip 6pt
{\bf Theorem 2.3.} {\it Let $\tau\in \r_-^{\times} $. Then (independent
of $\tau$) $N_{\tau} = R$ so that (2.6) becomes $$O_{\tau}\cong
N/R.\eqno (2.18)$$}
\indent {\bf Proof.} Let $\beta\in {\cal B}$. Then, by strong
orthogonality, 
$$[e_{\beta},\r_-]=
\Bbb C\,h_{\beta}\s \hh\eqno (2.19) $$ so that $\hbox{coad}\,\,
e_{\beta}(\r_-)= 0$. Hence $$\r\s \n_{\tau}.\eqno (2.20)$$ Conversely
let $v\in
\n_{\tau}$. We must show that $v\in\r$. Assume not. Then we may
assume that $v\in \ss$ and $v\neq 0$. Now by (2.10) we may write
$v =\sum_{\beta\in {\cal B}}\,v(\beta)$ where $v(\beta)\in\ss(\beta)$.
Let ${\cal B}_v =\{\beta \in {\cal B}\mid v(\beta)\neq 0\}$. Then ${\cal
B}_v$ is not empty. But by Lemma 2.2, $\hbox{coad}\,v(\beta)(\tau)\neq 0$ and
$\hbox{coad}\,v(\beta)(\tau)\in \ss_{-\beta}$ for 
$\beta\in {\cal B}_v$. Hence
$\hbox{coad}\,v(\tau)\neq 0$ by (2.11). This is a contradiction. \hfill QED

\vskip 6pt
{\bf Remark 2.4.} In effect Theorem 2.3.
depends on the fact, established in the proof, that if $0\neq
v\in \ss$ and $\tau \in \r_-^{\times} $, then $$0\neq \hbox{coad}\,v(\tau)\in
\ss_-.\eqno (2.21)$$
\vskip .1pc {\bf Theorem 2.5.} {\it If
$\tau,\tau'\in \r_-^{\times} $ are distinct,  then $O_{\tau}\cap O_{\tau'}=
\emptyset$ so that one has a disjoint union
$$\hbox{\rm Coad}\,N(\r_-^{\times} ) = \cup_{\tau\in \r_-^{\times} }\,O_{\tau}.\eqno (2.22)$$}
\indent {\bf
Proof.} It suffices to show that $$O_{\tau}\cap \r_-^{\times}  = \{\tau\}.\eqno
(2.23)$$ Let $u\in N$. We may write $u= \hbox{exp}\,z$ where $z\in \n$.
If
$z\in \r=\n_{\tau}$, then $\hbox{Coad}\,u(\tau) = \tau$. Hence we may
assume that $z\notin \r$ so that we can write $z = x + v$ where
$x\in \r$ and $0\neq v\in \ss$. By Remark 2.3, if
$y=\hbox{coad}\,z(\tau) $, then $0\neq y \in \ss_-$. Thus we may
write $$y = \sum_{\varphi\in
\Delta_+\setminus{\cal B}}\,c_{\varphi}\,e_{-\varphi}.\eqno (2.24)$$ Let
$\varphi'\in \Delta_+\setminus{\cal B}$ be a maximum element, relative to
the ordering (1.1), such that $c_{\varphi'}\neq 0$. But then it is
clear, from the exponentiation, that if we write 
$$\hbox{Coad}\,u(\tau)-\tau =
\sum_{\varphi\in
\Delta_+}\,d_{\varphi}\,e_{-\varphi}\eqno (2.25)$$ one has $d(\varphi') =
c(\varphi')$. Hence $\hbox{Coad}\,u(\tau)\notin \r_-$. In particular 
$Coad\,u(\tau)\notin T$. \hfill QED

\vskip 8pt {\bf 2.2.} Let $H\s G$ be
the subgroup corresponding to $\hh$ so that $B= N\,H$ is a Borel
subgroup of $G$ and $\b = \hbox{Lie}\,B$. Since $B$ normalizes $N$ the
dual space $\n_-$ to $\n$ is a $B$-module where $H$ operates via the
adjoint representation and of course $N$ operates via the
coadjoint representation. Obviously the decompositions (2.7) and
(2.8) are preserved by the action of $H$. Furthermore from the linear
independence of the elements of ${\cal B}$ one notes

 \vskip 6pt {\bf
Remark 2.6.} The Zariski open subvariety $\r_-^{\times} \s \r_-$ is stable under the action of $H$
and in fact $H$ operates transitively on $\r_-^{\times} $ so that $\r_-^{\times} $ is
isomorphic to a homogeneous space for $H$.

 \vskip 6pt The action
of $H$ on $\r_-^{\times} $ extends to an action of $H$ on the corresponding set
$\{O_{\tau},\,\tau\in \r_-^{\times} \}$ of $N$-coadjoint orbits. Since $H$
normalizes $N$ the following statement is obvious.

 \vskip 6pt {\bf
Proposition 2.7.} {\it For any
$\tau\in \r_-^{\times} $ and
$a\in H$ one has
$$O_{{\hbox{\sevenrm Ad}}\,a(\tau)} = \hbox{\rm Ad}\,a\,(O_{\tau}).\eqno (2.26)$$}
\indent Let
$$X=\cup_{\tau\in \r_-^{\times} }\,O_{\tau}\eqno (2.27)$$ so that, by 
Theorem 2.5,
the union (2.27) is disjoint.

\vskip 6pt Clearly $X$ (see Remark 2.6) is an
orbit of the action of $B$ on $\n_-$ so that $X$ has the structure
of an algebraic subvariety of $\n_-$ which is Zariski open in its
closure. On the other hand the product variety $N/R \times \r_-^{\times} $ is an
affine variety and one has a bijection $$\psi:N/R \times \r_-^{\times} \to
X,\eqno (2.28)$$ where if $[u]\in N/R$ denotes the left coset of
$u\in N$ in $N/R$, one has $$\psi(([u],\tau)) = \hbox{Coad}\,u(\tau).\eqno
(2.29)$$ But since $H$ normalizes both $N$ and $R$ it follows
easily that $N/R$ is a $B$-homogeneous space. The action of $H$ on
$\r_-^{\times} $ extends to an action of $B$ on $\r_-^{\times} $ where $N$ operates
trivially. Consequently $N/R \times \r_-^{\times} $ has the structure of a 
$B$-homogeneous space. 

\vskip 6pt {\bf Theorem 2.8.} {\it The map $\psi$
(see (2.28) and (2.29)) is a $B$-isomorphism of affine
$B$-homogeneous spaces. Furthermore $X$ is Zariski open in $\n_-$
so that $$\overline {X} = \n_-.\eqno (2.30)$$}
\vskip 3pt
\indent {\bf Proof.} We
have noted that $N/R \times \r_-^{\times} $ is an affine variety. Since a
homogeneous space of an affine algebraic group inherits a unique
algebraic structure, to establish the first statement, it suffices
only to see that
$\psi$ is a $B$-map. But this is immediate. But $dim\,\r_-^{\times}  = dim\,R$.
Thus $\hbox{dim}\,X =  \hbox{dim}\,\n_-$. This proves (2.30).\hfill QED

\vskip 8pt
\centerline{\bf 3. The characters of $H$ on $S(\n)^N$}\vskip 3pt
{\bf 3.1.} Let $\Lambda\s \hh^*$ be the weight lattice for
$(\hh,\g)$ and let $\Lambda_{ad}$ be the root sublattice of
$\Lambda$. For each $\nu\in \Lambda_{\hbox{\sevenrm{ad}}}$ let $\chi_{\nu}$ be the
character on $H$ defined so that for $a= \hbox{exp}\,x,\,x\in \hh$, one
has $\chi_{\nu}(a) = e^{\nu(x)}$. Of course the character group,
$\widehat {H}$, of $H$ is given by $$\widehat {H} = \{\chi_{\nu}\mid\nu
\in \Lambda_{\hbox{\sevenrm{ad}}}\}.$$ Recalling that $m=\hbox{card}\,{\cal B}$ let 
$\Lambda({\cal B}) \s \Lambda_{ad}$ be the free abelian group of
$\hbox{rank}\,m$, generated by ${\cal B}$, and let $$\widehat {H}({\cal B})
= \{\chi_{\nu}\mid \nu\in \Lambda({\cal B})\}.$$

If $V$ is a affine variety (over $\Bbb C$) we will let $A(V)$ denote
the affine algebra of regular functions on $V$. The quotient
field of $A(V)$, the algebra of rational functions on $V$, will be
denoted by $Q(V)$. If a linear algebraic group $G'$ operates 
algebraically on $V$, then $G'$ operates as a group of automorphisms
of
$A(V)$ so that if $g\in G',\, \phi\in A(V)$ and $v\in V$, then
$g\cdot\phi(v)= \phi(g^{-1}\cdot v)$. The group also operates as
a group of automorphisms of $Q(V)$ where, for $g\in
G',\,\phi,\phi'\in A(V)$ and $\phi'\neq 0$, then $g\cdot \phi/\phi'=
g\cdot\phi/g\cdot\phi'$. If $\g' = \hbox{Lie}\,G'$, then $\g'$ operates as a
Lie algebra of derivations of $A(V)$ and $Q(V)$. Using the fact
that, as one knows, $G'\cdot\phi$ spans a finite-dimensional
subspace of $A(V)$ for any $\phi\in A(V)$, one has $x\cdot\phi =
{d\over dt}(\hbox{exp}\,\,t\,x\cdot\phi)|_{t=0}$. If $0\neq \phi'\in A(V)$,
then
$x\cdot \phi/\phi' = (\phi'\,\,\,x\cdot \phi -\phi\,\,\,x\cdot
\phi')/ (\phi')^2$. If $M$ is any $G'$ module, then $M^{G'}$ will
the submodule of $G'$ invariants in $M$. 

Now the map $\psi$ (see (2.28)) induces a $B$-isomorphism
$$A(X) \to A(N/R)\otimes A(\r_-^{\times} )\eqno (3.1)$$ by Theorem 2.8. But then,
noting the action of $N$, the map (3.1) defines an $H$-isomorphism
$$A(X)^N\to A(\r_-^{\times} ).\eqno (3.2)$$ But then, recalling the affine
algebra of a (complex) torus one immediately has

\vskip 3pt {\bf Theorem 3.1.} {\it For any $\nu\in \Lambda({\cal B})$
there exists a unique (up to scalar multiplication) nonzero
$H$-weight vector $\xi_{\nu}$ in $A(X)^N$ of weight $\nu$. That
is, for any $a\in H$, $$a\cdot \xi_{\nu} =
\chi_{\nu}(a)\,\xi_{\nu}.\eqno (3.3)$$ Moreover the set,
$\{\xi_{\nu}\mid \nu\in \Lambda({\cal B})\}$, of weight vectors
are a basis of $A(X)^N$. In particular every weight in $A(X)^N$
occurs with multiplicity one and only weights in $\Lambda({\cal
B})$ occur.}

\vskip 6pt {\bf 3.2.} The Killing form pairing of 
$\n$ and $\n_-$ identifies the symmetric algebra $S(\n)$ with
$A(\n_-)$ and the quotient field $F(\n)$ of $S(\n)$ with
$Q(\n_-)$. Of course $S(\n)$ is a unique factorization domain. 

\vskip 3pt
{\bf Remark 3.2.} It is evident that if $\phi\in S(\n)$ is a
prime polynomial and $b\in B$, then $b\cdot \phi$ is again a
prime polynomial and $\psi = \phi_1\cdots \phi_k$ is the prime
factorization of $0\neq \psi\in S(\n)$, then $$ b\cdot \psi = (b\cdot
\phi_1)\cdots (b\cdot \phi_k)\eqno (3.4)$$ is the prime
factorization of $b\cdot \psi$.

\vskip3pt {\bf Proposition 3.3.} {\it
$Q(\n_-)^N$ is the quotient field of
$S(\n)^N$. Furthermore the prime factors of any $0\neq \psi\in
S(\n)^N$ are also in $S(\n)^N$.} \vs {\bf Proof.} Let $0\neq
\gamma\in Q(\n_-)^N$. Write $\gamma = \phi/\psi$ where
$\phi,\psi\in S(\n)$. Then $\psi\,\gamma = \phi$ and hence for any
$u\in N$ one has $$(u\cdot \psi)\,\gamma = u\cdot \phi.\eqno (3.5)$$
But $N\cdot \psi$ spans a finite-dimensional $N$-submodule $M$ of 
$S(\n)$. By the unipotence of $N$ and its action on $M$ there
exists $0\neq	 \psi'\in M^N\s S(\n)^N$. But then if $\phi' =
\psi'\,\gamma$, it follows from (3.5) that $\phi'\in S(\n)^N$. But 
$\gamma = \phi'/\psi'$. This proves the first statement of the
proposition. 

Now let $0\neq \psi\in S(\n)^N$ and let $\psi = \phi_1\cdots\phi_k$
be a prime factor decomposition of $\psi$. But then for any $u\in
N$, $$\psi = (u\cdot \phi_1)\cdots (u\cdot \phi_k)\eqno (3.6)$$ is
another prime factor decomposition of $\psi$. By the continuity of
the action of $N$ and the uniqueness (up to scalar multiplication)
of the prime factor decomposition, for any $j\in \{1,\ldots,k\}$,
there exists $c\in \Bbb C^{\times}$ such that $u\cdot \phi_j =
c\,\phi_j$. But, by unipotence, $1$ is the only eigenvalue of the
action of $u$ on $S(\n)$. Thus $c=1$ and hence $\phi_j\in S(\n)^N$.
\hfill QED

\vs Now, by (2.30), one has $Q(X) = Q(\n_-)$. Thus we have the
$B$-inclusions $$S(\n)\s A(X)\s Q(\n_-).\eqno (3.7)$$

\vskip .2pc {\bf
Remark 3.4.} Note that in (3.7), as functions on $\n_-$, the
elements of $S(\n)$ are exactly the functions in $A(X)$ which
extend, as regular functions, (i.e., everywhere defined rational
functions) to all of
$\n_-$. 

\vskip 6pt But now (3.7) yields the inclusion
$$S(\n)^N \s A(X)^N\eqno (3.8)$$ of (completely reducible)
$H$-modules. Recalling Theorem 3.1 let $$\Lambda_{\n}({\cal B})
= \{\nu\in \Lambda({\cal B})\mid \xi_{\nu}\in S(\n)^N\}.\eqno (3.9)$$
Also let $\Lambda_{\hbox{\sevenrm{dom}}}({\cal B})$ be the set of all dominant
weights in $\Lambda({\cal B})$. In the following result it is proved that $\Lambda_{\n}({\cal
B})\s \Lambda_{\hbox{\sevenrm{dom}}}({\cal B})$.  It is established as argued in the Introduction that in
fact 
$\Lambda_{\n}({\cal
B})= \Lambda_{\hbox{\sevenrm{dom}}}({\cal
B})$. 

\vskip 6pt{\bf Theorem 3.5.} {\it Every
$H$-weight in $S(\n)^N$ occurs with multiplicity 1. Moreover
$\Lambda_{\n}({\cal B})$ is the set of such $H$-weights. 
Furthermore $$ \{\xi_{\nu}\mid \nu\in \Lambda_{\n}({\cal B})\}\eqno
(3.10)$$ is a basis of $S(\n)^N$. Finally $$\Lambda_{\n}({\cal
B})\s \Lambda_{\hbox{\sevenrm{dom}}}({\cal B}).\eqno (3.11)$$} 
\indent {\bf Proof.}
Except for (3.10) the theorem follows immediately from (3.8) and
Theorem 3.1. The inclusion (3.11) follows immediately from the
fact that $\xi_{\nu}$ for $\nu\in \Lambda_{\n}({\cal B})$ is
necessarily the highest weight vector for the $\g$-submodule, in the
symmetric algebra $S(\g)$, generated by $\xi_{\nu}$. \hfill QED

\vskip 6pt Let
$${\cal P} = \{ \nu \in \Lambda_{\n}({\cal B})\mid \hbox {$\xi_{\nu}$ is a prime polynomial in
$S(\n)$}\}.\eqno (3.12)$$

\vskip .5pc {\bf Theorem 3.6.} {\it One has
$\hbox{card}\,{\cal P} = m$ where, we recall $m=\hbox{card}\,{\cal B}$, so that
we can write $${\cal P} = \{\mu_1,\ldots, \mu_{m}\}.\eqno (3.13)$$
Furthermore the weights $\mu_i$ in ${\cal P}$ are
linearly independent and the set $P$ of prime polynomials,
$\xi_{\mu_i},\,i=1,\ldots, m$, are algebraically independent. In
addition one has a bijection $$\Lambda_{\n}({\cal B}) \to (\Bbb
N)^m,\quad\nu \mapsto (d_1(\nu),\ldots, d_m(\nu))\eqno (3.14)$$
such that, writing $d_i= d_i(\nu)$, up to scalar multiplication,
$$\xi_{\nu} = \xi_{\mu_1}^{d_1}\cdots \xi_{\mu_m}^{d_m}\eqno
(3.15)$$ and (3.15) is the prime factorization of $\xi_{\nu}$ for
any $\nu \in \Lambda_{\n}{\cal B}$. Finally $$S(\n)^{N} = \Bbb
C[\xi_{\mu_1},\ldots,\xi_{\mu_m}]\eqno (3.16)$$ so that $S(\n)^{N}$
is a polynomial ring in $m$ generators.}

\vskip 6pt {\bf Proof.} Let $\nu
\in \Lambda_{\n}({\cal B}) $ and consider the prime factorization of
$\xi_{\nu}$. We use the notation of Remark 3.2 where 
$\psi =\xi_{\nu} $ and $b\in H$. Then since $b\cdot \xi_{\nu} =
\chi_{\nu}(b)\,\xi_{\nu}$ it follows from (3.4) that the right
side of (3.4) is another prime factorization of $\xi_{\nu}$. By the
continuity of the action of $H$ and the uniqueness of
the factorization it follows that for any $j=1,\ldots,k$, there
exists $\chi_{j}\in {\widehat H}$ such that $b\cdot \phi_j=
\chi_j(b)\phi_j$ for all $b\in H$. But then, by Theorem 3.1
and Theorem 3.5, one has $\chi_j = \chi_{\nu_j}$ for a unique
$\nu_j\in \Lambda_{\n}(B)$ and up to scalar multiplication
$\phi_j = \xi_{\nu_j}$. But also $\nu_j\in {\cal P}$. Thus for
one thing this shows ${\cal P}$ is not empty and up to scalar
multiplication
$$\xi_{\nu} = \xi_{\nu_1}\cdots \xi_{\nu_k}\eqno (3.17)$$ and
(3.17) is the prime factorization of $\xi_{\nu}$.

For a positive integer $r\leq \hbox{card}\,{\cal P}$ let
$\mu_1',\ldots,\mu_{r}'$ be
$r$ distinct elements of ${\cal P}$. For any
$d=(d_1,\ldots,d_{r})\in (\Bbb Z_+)^{r}$ let
$\nu(d)\in
\Lambda({\cal B})$ be defined by putting $\nu(d)= \sum_{j=1}^{r}\,
d_j\,\mu_j'$. But, up to scalar multiplication,
$$\xi_{\nu(d)} = \xi_{\mu_1'}^{d_1}\cdots
\xi_{\mu_{r}'}^{d_r}\eqno(3.18)$$ so that $\nu(d)\in 
\Lambda_{\n}({\cal B})$ and (3.18) is the prime factorization of 
$\xi_{\nu(d)}$. One also has a map
$$(\Bbb Z_+)^{r}\to
\Lambda_{\n}({\cal B}),\qquad d \mapsto \nu(d).\eqno (3.19)$$ But
then, by the uniquess of the prime factorization, the map (3.19)
is necessarily injective. But weight vectors belonging to distinct
weights are linearly independent (see also Theorem 3.5) so that
 no nontrivial linear combination of the monomials on the right
side of (3.18) can vanish. This proves $$\xi_{\mu_1'},\ldots,
\xi_{\mu_r'}\,\,\hbox{are algebraically independent}.\eqno (3.20)$$
But by (2.5) and (3.2) the transcendence degree of $A(X)^N$ is $m$.
But $$S(\n)^N\s A(X)^N\s Q(\n_-)^N\eqno (3.21)$$ by (3.6). Hence
$$\hbox{the transcendence degree of $S(\n)^N$ is $m$}\eqno (3.22)$$
by Proposition 3.2. Thus
$r\leq m$ and hence if $n = \hbox{card}\,{\cal P}$ one has $n\leq m$. But
then if we choose
$r = n$, the map (3.19) is surjective by (3.17) and one must have
$S(\n)^N = \Bbb C[\xi_{\mu_1'},\ldots,\xi_{\mu_n'}]$. But then
$n=m$ by (3.21). Except for the statement that $\nu_i$ are linearly
independent weights this proves Theorem 3.6, noting that (3.14) is
the inverse of the bijection (3.19) when we choose $\mu_j= \mu_j'$.
Assume that the
$\mu_i$ are not linearly dependent. Since these weights lie in a
lattice there exists, over the rational numbers, a vanishing
nontrivial linear combination. Clearing
denominators such a linear combination exists over the integers. But
this implies that there exists distinct $d,d'\in (\Bbb N)^m$ such
that $\nu(d) = \nu(d')$. This contradicts the bijectivity of
(3.14).\hfill QED

\vskip 6pt {\bf Remark 3.7.} We note the following uniqueness
statement. If
$\{\nu_1,\ldots ,\nu_k\}$ is any subset of $\Lambda_{\n}({\cal B})$
with the property that there exists a bijection
$$\Lambda_{\n}({\cal B})\to ({\Bbb Z_+})^k,\quad \nu 
\mapsto (e_1(\nu),\ldots,e_k(\nu))$$ such that,
writing $e_i=e_i(\nu)$, $$\xi_{\nu} = \xi_{\nu_1}^{e_1}\cdots
\xi_{\nu_k}^{e_k}$$ up to scalar multiplication, then one must have
$k = m$ and
$\{\nu_1,\ldots ,\nu_k\}$ is some reordering of ${\cal P} =
\{\mu_1,\ldots,\mu_m\}$. This is clear since otherwise one would
have a contradiction of the primeness of the $\xi_{\mu_j}$. In
particular the generators constructed by A. Joseph in [J] of
(without any reference to primeness) must necessarily be the same as
our $\xi_{\nu_j}$.

\vskip 6pt  As to all the weights $\nu$ in the group
(not semigroup) 
$\Lambda ({\cal B})$ (see Theorem 3.1), one has an extension of the map
(3.14) involving $\Bbb Z$ (and hence negative integers as well). 

\vskip 6pt
{\bf Theorem 3.8.} {\it There exists a bijection
$$\Lambda({\cal B})\to ({\Bbb Z})^m,\quad \nu\mapsto
(d_1(\nu),\ldots, d_{m}(\nu))\eqno (3.23)$$ such that, writing $d_i
= d_i(\nu)$, $$\xi_{\nu} = \xi_{\mu_1}^{d_1}\cdots
\xi_{\mu_k}^{d_k}\eqno (3.24)$$ up to scalar multiplication,
recalling that the
$\xi_{\nu}$, for
$\nu\in \Lambda({\cal B})$ is a basis of $A(X)^N$. See Theorem
3.1.}

\indent{\bf Proof.} Let $\nu \in \Lambda({\cal B})$. By (3.21) and
Proposition 3.3, up to scalar multiplication, one can uniquely write
$\xi_{\nu} = p/q$ where $p,q \in S(\n)^N$ and $p$ and $q$ are prime
to one another. But then if $a\in H$ one has $\chi_{\nu}(a) = a\cdot
p/a\cdot q$. By uniqueness both $p$ and $q$ must be weight vectors
in $S(\n)^N$. But then Theorem 3.8 follows from Theorem 3.6.\hfill QED

\vskip 6pt
{\bf Remark 3.9.} It follows from Theorem 3.8 that both ${\cal B}$
and ${\cal P}$ are bases for the free abelian group
$\Lambda({\cal B})$. Consequently there must be a matrix in
$Sl(m,\Bbb Z)$ which expresses one such basis in terms of the
other.

 \vskip 6pt
For $\nu\in \hh^*$ and $\beta\in {\cal B}$ let $r_{\beta}(\nu) = (\beta,\nu)/(\beta,\beta)$ so that,
by (1.7), $$\eqalign{w_o(\nu) = -\nu\,\,&\iff\nu\,\,\hbox{is in the span of ${\cal B}$}\cr &\iff \nu
= 
\sum_{\beta\in {\cal B}} r_{\beta}(\nu)
\,\beta.\cr}\eqno (3.25)$$  
In any case put 
$$r(\nu) = \sum_{\beta\in {\cal B}} r_{\beta}(\nu).\eqno (3.26)$$ Let $\Lambda_{dom}$ be the set
of all dominant weights so that

 \vskip .5pc {\bf Remark 3.10.} If
$\nu\in \Lambda_{\hbox{\sevenrm{dom}}}$ note that, for all $\beta\in {\cal B}$, $r_{\beta}(\nu)\in \Bbb Z_+/2$
and if
$\nu\in \Lambda_{\hbox{\sevenrm{dom}}}({\cal B})$ then $w_o\nu = -\nu$ and $r_{\beta}(\nu)\in \Bbb Z_+$. In
particular this is true for
$\nu\in \Lambda_{\n}({\cal B})$ by (3.11).
 
\vskip 8pt {\bf 3.3.} Let
$\Gamma$ be the set of all maps
$\gamma:\Delta_+\to \Bbb Z_+$. For $\gamma \in \Gamma$ let $d(\gamma) =
\sum_{\varphi\in \Delta_+} \gamma(\varphi)$ and let $\gamma^{\hbox{\sevenrm{root}}}$ be the
element in the root lattice given by putting $\gamma^{\hbox{\sevenrm{root}}} = \sum_{\varphi\in
\Delta_+} \gamma(\varphi)\,\varphi$. Also let $z_{\gamma}\in S^{d(\gamma)}(\n)$
be given by putting $$z_{\gamma} = \prod_{\varphi\in \Delta_+}
e_{\varphi}^{\gamma(\varphi)}.\eqno (3.27)$$ Let $\nu \in \Lambda_{\n}({\cal
B})$. Now from the multiplicity 1 condition (see Theorem 3.1) it follows that
$$\xi_{\nu},\,\,\hbox{for any
$\nu \in
\Lambda_{\n}({\cal B})$, is a homogeneous polynomial in $S(\n)$}.\eqno (3.28)$$ 
Let $\hbox{deg}\,\nu$ be the degree of the homogeneous polynomial $\xi_{\nu}$ and let 
$$\Gamma(\nu) =
\{\gamma\in
\Gamma\mid d(\gamma) =\hbox{deg}\,\nu\,\,\hbox{and}\,\,\gamma^{\hbox{\sevenrm{root}}} = \nu\}.$$ But then 
there exists scalars $s_{\gamma},\,\gamma\in \Gamma(\nu)$ such that
$$\xi_{\nu}= \sum_{\gamma \in \Gamma(\nu)}\,s_{\gamma}\,z_{\gamma}.\eqno
(3.29)$$
\vskip .5pc But now since $X$ is Zariski open in $\n_-$ (see Theorem
2.8) there exists $x\in X$ such that $\xi_{\nu}(x)\neq 0$. But $x=
\hbox{Coad}\,b^{-1}(t)$ for $b\in N$ and $t\in \r_-^{\times} $. Write $$t = \sum_{\beta\in {\cal
B}}t_{\beta}\,e_{-\beta},\eqno (3.30)$$ where $t_{\beta}\in \Bbb C^{\times}$ for
all $\beta\in {\cal B}$. But by $N$-invariance $$\eqalign{0&\neq
\xi_{\nu}(x)\cr&=\xi_{\nu}(t).\cr}\eqno (3.31)$$ But $z_{\gamma}(t) = 0$ for $\gamma\in
\Gamma(\nu)$ unless $\gamma = \gamma_{\r}$ where $\gamma_{\r}(\varphi)=0$ for
$\varphi\notin {\cal B}$. But any such $\gamma_{\r}$ is clearly unique where
necessarily
$\gamma(\beta) = r_{\beta}(\nu)$, for all $\beta\in {\cal B}$. Hence, also one
must have
$deg\,\nu = r(\nu)$.  We have in fact thus proved 

\vskip 6pt {\bf Theorem 3.11.} {\it Let
$\nu\in \Lambda_{\n}({\cal B})$. Then $\hbox{deg}\,\nu=r(\nu)$. Furthermore (see (3.29)) one must have
$s_{\gamma_{\r}}\neq 0$ and if $t\in \r_-^{\times} $ is given by (3.30), then $$\xi_{\nu}(t) =
s_{\gamma_{\r}}\prod_{\beta \in {\cal B}}t_{\beta}^{r_{\beta}(\nu)}.\eqno
(3.32)$$}
\vskip 1 pc
\centerline{\bf References}\vskip 6pt
\item {[D]} Jacques Dixmier, {\it Sur les repr\'esentations unitaires
des groupes de Lie nilpotent.} IV, Canad. J. Math.,{\bf 11}(1959),
321--344.
\item {[LW]} Ronald Lipsman and Joseph Wolf, {\it Canonical
semi-invariants and the Plancherel formula for parabolic groups},
Trans. Amer. Math. Soc., {\bf 269}(1982), 111--131.
\item {[J]} Anthony Joseph, {\it A preparation theorem for the
prime spectrum of a semisimple Lie algebra}, Journ. of Alg., {\bf
48}, No.2, (1977), 241--289.
\item {[K]} Bertram Kostant, {\it $\hbox{Cent}\, U(\n)$ and a construction
of Lipsman-Wolf}, arXiv /0178653, [math.RT], Jan 12, 2011.
\vskip 1pc
\noindent Bertram Kostant (Emeritus)

\noindent Department of Mathematics

\noindent Cambridge, MA 02139

\noindent kostant@math.mit.edu

\bye
\end

\end